\documentclass[12pt, reqno]{amsart}
\usepackage{amsmath, amsthm, amscd, amsfonts, amssymb, graphicx, color}
\usepackage[bookmarksnumbered, colorlinks, plainpages]{hyperref}

\newtheorem{theorem}{Theorem}[section]
\newtheorem{lemma}[theorem]{Lemma}

\newtheorem{corollary}[theorem]{Corollary}
\theoremstyle{definition}
\newtheorem{definition}[theorem]{Definition}
\newtheorem{example}[theorem]{Example}

\theoremstyle{remark}
\newtheorem{remark}[theorem]{Remark}
\numberwithin{equation}{section}

\allowdisplaybreaks
\begin{document}

\title[Continuous frames in Hilbert $C^*$-modules]
{Continuous frames in Hilbert $C^*$-modules}

\author[H.Ghasemi]{Hadi Ghasemi }
\address{Hadi Ghasemi \\ Department of Mathematics and Computer
Sciences, Hakim Sabzevari University, Sabzevar, P.O. Box 397, IRAN}
\email{ \rm h.ghasemi@hsu.ac.ir}
\thanks{h.ghasemi@hsu.ac.ir (Hadi Ghasemi)}
\author[T.L. Shateri]{Tayebe Lal Shateri }
\address{Tayebe Lal Shateri \\ Department of Mathematics and Computer
Sciences, Hakim Sabzevari University, Sabzevar, P.O. Box 397, IRAN}
\email{ \rm  t.shateri@hsu.ac.ir; shateri@ualberta.ca}
\thanks{*The corresponding author:
t.shateri@hsu.ac.ir; shateri@ualberta.ca (Tayebe Lal Shateri)}
 \date{07, 04, 2022}
 \maketitle

\begin{abstract}
In the present paper the notion of continuous frames is introduced and some results of these frames are proved. Next, we give the concept of duals of continuous frames in Hilbert $C^*$-modules and investigate some properties of them.
\vskip 3mm
\end{abstract}

\section{Introduction And Preliminaries}
Frames for Hilbert spaces were introduced in 1952 by Duffin and Schaeffer \cite{DS} to study some problems in nonharmonic Fourier series and widely studied from 1986 since the great work by Daubechies, Grossmann and Meyer \cite{DG}. Now, frames have been widely applied in signal processing, sampling, filter bank theory, system modeling, Quantum information, cryptography, etc. \cite{BL,EL,FE,ST}.  Various generalizations of frames e.g. frames of subspaces, wavelet frames, $g$-frames, weighted and controlled frames were developed, see  \cite{CA,SU,BAL}. We refer to \cite{CH} for an introduction to frame theory and its applications. The concept of a generalization of frames to a family indexed by some locally compact space endowed with a Radon measure was proposed by G. Kaiser \cite{KAI} and independently by Ali, Antoine and Gazeau \cite{AAG}. These frames are known as continuous frames.

Frank and Larson \cite{FL} extended the notion of a frame for an operator on a Hilbert  $C^*$-module. For a discussion of frames in Hilbert $C^*$-modules, we refer to Refs. \cite{AR,GHSH1,JI,SH}. Hilbert $C^*$-module is a generalization of a Hilbert space that allows the inner product to take values in a $C^*$-algebra rather than the field of complex numbers.
The extended results of this more general framework are not a routine generalization, because there are essential differences between Hilbert $C^*$-module and Hilbert Space.  For example, we know that every bounded operator on a Hilbert space has an unique adjoint, while this fact does not hold for bounded operators on a Hilbert $C^*$-module. Also,  any closed subspace in a Hilbert space has an orthogonal complement, but it is not true, in general, for a Hilbert $C^*$-module. 

The paper is organized as follows. First, we recall the basic definitions and some notations about Hilbert $C^*$-modules, and we also give some properties of them which we will use in the later sections. In Section 2, we introduce the notion of continuous frames in Hilbert $C^*$-modules. We present some results of frames in the view of continuous frames. In section 3, we define duals of continuous frames in Hilbert $C^*$-modules and investigate some properties of them.

First, we recall some definitions and basic properties of Hilbert $C^*$-modules. We give only a brief introduction to the theory of Hilbert $C^*$-modules to make our explanations self-contained. For comprehensive accounts we refer to  \cite{LAN,MT,OLS}.\\
Throughout this paper, $\mathcal A$ is a unital $C^*$-algebra.
\begin{definition}
A \textit{pre-Hilbert module} over a unital $C^*$-algebra $\mathcal A$ is a complex vector space $U$ which is also a left $\mathcal A$-module equipped with an $\mathcal A$-valued inner product $\langle .,.\rangle :U\times U\to \mathcal A$ which is $\mathbb C$-linear and $\mathcal A$-linear in its first variable and satisfies the following conditions:\\
$(i)\; \langle f,f\rangle \geq 0$,\\
$(ii)\; \langle f,f\rangle =0$  iff $f=0$,\\
$(iii)\; \langle f,g\rangle ^*=\langle g,f\rangle ,$\\
$(iv)\; \langle af,g\rangle=a\langle f,g\rangle ,$\\
for all $f,g\in U$ and $a\in\mathcal A$.
\end{definition}
A pre-Hilbert $\mathcal A$-module $U$ is called \textit{Hilbert $\mathcal A$-module} if $U$ is complete with respect to the topology determined by the norm $\|f\|=\|\langle f,f\rangle \|^{\frac{1}{2}}$.

By \cite[Example 2.46]{JI}, if $\mathcal A$ is a $C^*$-algebra, then it is a Hilbert $\mathcal A$-module with respect to the inner product
$$\langle a,b\rangle =ab^*,\quad (a,b\in \mathcal A).$$
\begin{example}
Let $l^2(\mathcal A)$ be the set of all sequences $\{a_n\}_{n\in \mathbb N}$ of elements of a $C^*$-algebra $\mathcal A$ such that the series $\sum_{n=1}^{\infty}a_na_n^*$ is convergent in $\mathcal A$. Then $l^2(\mathcal A)$ is a Hilbert $\mathcal A$-module with respect to the pointwise operations and inner product defined by
\begin{equation*}
\langle \{a_n\}_{n\in \mathbb N},\{b_n\}_{n\in \mathbb N}\rangle =\sum_{n=1}^{\infty}a_nb_n^*.
\end{equation*}
\end{example}
In the following lemma the \textit{Cauchy-Schwartz inequality} extends to Hilbert $C^*$-modules.
\begin{lemma}
\cite[Lemma 15.1.3]{OLS} (\textbf{Cauchy-Schwartz inequality}) Let $U$ be a Hilbert $C^*$-modules over a unital $C^*$-algebra $\mathcal A$. Then
\begin{equation*}
\|\langle f,g\rangle \|^{2}\leq\|\langle f,f\rangle \|\;\|\langle g,g\rangle \|
\end{equation*}
for all $f,g\in U$.
\end{lemma}

\begin{definition}
Let $U$ and $V$ be two Hilbert $C^*$-modules over a unital $C^*$-algebra $\mathcal A$. A map $T:U\to V$ is said to be \textit{adjointable} if there exists a map $T^{*}:V\to U$ satisfying
$$\langle Tf,g\rangle =\langle f,T^*g\rangle $$
for all $f\in U, g\in V$. Such a map $T^*$ is called the \textit{adjoint} of $T$. By $End_{\mathcal A}^*(U)$ we denote the set of all adjointable maps on $U$.
\end{definition}
It is surprising that every adjointable operator is automatically linear and bounded.\\
We need the following results in next sections.
\begin{lemma} \label{UKR}
\cite[Lemma 1.1]{XS} Let $U$ and $V$ be two Hilbert $C^*$-modules over a unital $C^*$-algebra $\mathcal A$ and $T\in End_{\mathcal A}^*(U,V)$ have closed range. Then $T^*$ has closed range and 
$$U=Ker(T)\oplus R(T^{*})\;\;\;\;\; ,\;\;\;\;\; V=Ker(T^{*})\oplus R(T)$$
\end{lemma}
\begin{lemma}\label{SB}
\cite{AR} Let $U$ and $V$ be two Hilbert $C^*$-modules over a unital $C^*$-algebra $\mathcal A$ and $T\in End_{\mathcal A}^*(U,V)$. Then the following are equivalent:\\
(i)\;$T$ is surjective\\
(ii)\;$T^{*}$ is bounded below with respect to the norm; i.e.
\begin{equation*}
\exists m>0\;\;\;\;\;s.t\;\;\;\;\;\|T^{*}f\|\geq m\| f\|.
\end{equation*}
(iii)\;$T^{*}$ is bounded below with respect to the inner product; i.e.
\begin{equation*}
\exists m>0\;\;\;\;\;s.t\;\;\;\;\;\langle T^{*}f,T^{*}f\rangle\geq m\langle f,f\rangle.
\end{equation*}
for all $f\in V$.
\end{lemma}
\begin{theorem}\label{BL}
\cite[Theorem 2.1.4]{MT} Let $U$ and $V$ be two Hilbert $C^*$-modules over a unital $C^*$-algebra $\mathcal A$ and $T\in End_{\mathcal A}^*(U,V)$. Then The following are equivalent:\\
(i)\;$T$ is bounded and $\mathcal A$-linear.\\
(ii)\;There exists $k>0$ such that
\begin{equation*}
\langle Tf,Tf\rangle\leq k\langle f,f\rangle.
\end{equation*}
for all $f\in U$.
\end{theorem}
\section{Continuous frames in Hilbert $C^*$-modules}
In this section, we introduce continuous frames in Hilbert $C^*$-modules over a unital $C^*$-algebra $\mathcal A$, and then we give some results for these frames. 

Let $\mathcal Y$ be a Banach space, $(\mathcal X,\mu)$ a measure space, and $f:\mathcal X\to \mathcal Y$ a measurable
function. The integral of the Banach-valued function $f$ has been defined by Bochner and others. Most properties of this integral are similar to those of the integral of real-valued functions (see \cite{DAN,YOS}). Since every $C^*$-algebra and Hilbert $C^*$-module is a Banach space, hence we can use this integral in these spaces.

In the following, we assume that $\mathcal A$ is a unital $C^*$-algebra, $U$ is a Hilbert $C^*$-module over $\mathcal A$ and $(\Omega ,\mu)$ is a measure space.
\begin{definition}
Let $(\Omega ,\mu)$ be a measure space and $\mathcal A$ is a unital $C^*$-algebra. We define,
\begin{equation*}
L^{2}(\Omega ,A)=\lbrace\varphi :\Omega \to A\quad ;\quad \Vert\int_{\Omega}\vert(\varphi(\omega))^{*}\vert^{2} d\mu(\omega)\Vert<\infty\rbrace 
\end{equation*}
For any $\varphi ,\psi \in L^{2}(\Omega , A)$, the inner product is defined by $\langle \varphi ,\psi\rangle = \int_{\Omega}\langle\varphi(\omega),\psi(\omega)\rangle d\mu(\omega)$ and the norm is defined by $\|\varphi\|=\|\langle \varphi,\varphi\rangle \|^{\frac{1}{2}}$. It was shown in \cite{LAN} that $L^{2}(\Omega , A)$ is a Hilbert $\mathcal A$-module.
\end{definition}
Now we define continuous frames in Hilbert $\mathcal A$-modules.
\begin{definition}
A mapping $F:\Omega \to U$ is called a continuous frame for $U$ if\\
$(i)\; F$ is weakly-measurable, i.e, for any $f\in U$, the mapping $\omega\longmapsto\langle f,F(\omega)\rangle $ is measurable on $\Omega$.\\
$(ii)$ There exist constants $A,B>0$ such that
\begin{equation}\label{eq1}
A\langle f,f\rangle \leq \int_{\Omega}\langle f,F(\omega)\rangle \langle F(\omega),f\rangle d\mu(\omega)\leq B\langle f,f\rangle  ,\quad (f\in U). 
\end{equation}
\end{definition}

The constants $A,B$ are called \textit{lower} and \textit{upper} frame bounds, respectively. The mapping $F$ is called \textit{Bessel} if the right inequality in \eqref{eq1} holds and is called \textit{tight} if  $A=B$.
\begin{definition}
A continuous frame $F:\Omega \to U$ is called \textit{exact} if for every measurable subset $\Omega_{1}\subseteq\Omega$ with
$0<\mu(\Omega_{1})<\infty$, the mapping $F:\Omega\backslash\Omega_{1} \to U$ is not a continuous frame for $U$.
\end{definition}
Now, we give an example of a continuous frame in a Hilbert $C^*$-module.
\begin{example}\label{ex1}
Assume that $\mathcal A=\Big\{\begin{pmatrix}
a&0\\
0&b
\end{pmatrix}: x,y\in \mathbb C\Big\}$, then $\mathcal A$ is a unital $C^*$-algebra. Also $\mathcal A$ is a Hilbert $C^*$-module over itself, with the following inner product,
\begin{equation*}
\begin{array}{ll}
\langle .,.\rangle:\mathcal A\times \mathcal A\;\to \quad \mathcal A \\
\qquad\; (M,N)\longmapsto M(\overline{N})^t.
\end{array}
\end{equation*}
Suppose that $(\Omega ,\mu)$ is a measure space where  $\Omega=[0,1]$ and $\mu$ is the Lebesgue measure.  Consider the mapping $F:\Omega\to \mathcal A$ defined as  
$F(\omega)=\begin{pmatrix}
2\omega&0\\
0&\omega -1
\end{pmatrix}$, for any $\omega\in \Omega$.\\ 
For each $f=
\begin{pmatrix}
a&0\\
0&b
\end{pmatrix}
\in \mathcal A$, we have
\begin{align*}
\int_{\Omega}\langle f,F(\omega)\rangle\langle F(\omega),f\rangle d\mu(\omega)&=\int_{[0,1]}\langle\begin{pmatrix}
a&0\\
0&b
\end{pmatrix},\begin{pmatrix}
2\omega&0\\
0&\omega -1
\end{pmatrix}\rangle\langle\begin{pmatrix}
2\omega&0\\
0&\omega -1
\end{pmatrix} ,\begin{pmatrix}
a&0\\
0&b
\end{pmatrix}\rangle  d\mu (\omega)\\
&=\int_{[0,1]}
\begin{pmatrix}
2\omega a&0\\
0&(\omega -1)b
\end{pmatrix}\begin{pmatrix}
2\omega\overline{a}&0\\
0&(\omega -1)\overline{b}
\end{pmatrix} d\mu(\omega)\\
&=\int_{[0,1]}
\begin{pmatrix}
4\omega ^{2} &0\\
0&(\omega -1)^{2}
\end{pmatrix}\begin{pmatrix}
\vert a\vert ^2&0\\
0&\vert b\vert ^2
\end{pmatrix}d\mu(\omega)\\
&=\begin{pmatrix}
\vert a\vert ^2&0\\
0&\vert b\vert ^2
\end{pmatrix}\int_{[0,1]}
\begin{pmatrix}
4\omega ^{2} &0\\
0&(\omega -1)^{2}
\end{pmatrix}d\mu(\omega)=\begin{pmatrix}
\dfrac{4}{3}&0\\
0&\dfrac{1}{3}
\end{pmatrix}\begin{pmatrix}
\vert a\vert ^2&0\\
0&\vert b\vert ^2
\end{pmatrix},
\end{align*}
hence
\begin{equation*}
\dfrac{1}{3}\langle f,f\rangle\leq\int_{\Omega}\langle f,F(\omega)\rangle\langle F(\omega) ,f\rangle d\mu(\omega)\leq \dfrac{4}{3}\langle f,f\rangle.
\end{equation*}
Therefore, $F$ is a continuous frame with bounds $A=\dfrac{1}{3}$ and $B=\dfrac{4}{3}$.
\end{example}
Similar to continuous frames in Hilbert spaces, we introduce a pre-frame operator and frame operator for continuous frames in Hilbert $C^{\ast}$-modules.
\begin{definition}
Let $F:\Omega \to U$ be a continuous frame. Then \\
(i)\;The \textit{synthesis operator} or \textit{pre-frame operator} $T_{F}:L^{2}(\Omega , A)\;\to U$ weakly defined by
\begin{equation}
\langle T_{F}\varphi ,f\rangle =\int_{\Omega}\varphi(\omega)\langle F(\omega),f\rangle d\mu(\omega),\quad (f\in U).
\end{equation}
(ii)\; The adjoint of $T$, called The \textit{analysis operator} $T^{\ast}_{F}:U\;\to L^{2}(\Omega , A)$ is defined by
\begin{equation}
(T^{\ast}_{F}f)(\omega)=\langle f ,F(\omega)\rangle\quad (\omega\in \Omega).
\end{equation}
\end{definition}
 \begin{definition}
Let $F:\Omega \to U$ be a continuous frame for Hilbert $C^{\ast}$-module $U$. Then the frame operator $S_{F}:U\;\to U$ is weakly defined by
\begin{equation}
\langle S_{F}f ,f\rangle =\int_{\Omega}\langle f,F(\omega)\rangle\langle F(\omega),f\rangle d\mu(\omega),\quad (f\in U).
\end{equation}
\end{definition}
In following theorem, we investigate some properties of frame operators. 
\begin{theorem}\label{TH1}
Let $F:\Omega \to U$ be a continuous frame for Hilbert $C^{\ast}$-module $U$ with bounds $A,B$. Then the pre-frame operator $T_{F}:L^{2}(\Omega , A)\;\to U$ is well defined, surjective, adjointable $\mathcal A$-linear map and bounded with $\| T\|\leq\sqrt{B}$ . Moreover the analysis operator $T^{\ast}_{F}:U\;\to L^{2}(\Omega , A)$ is injective and has closed range.
\end{theorem}

\begin{proof}
Let $F:\Omega \to U$ be a continuous frame for Hilbert $C^{\ast}$-module $U$ with bounds $A,B$. Then 

(i)\; $T$ is adjointable and $T^{\ast}$ is its adjoint. Because for $f\in U$ and $\varphi\in L^{2}(\Omega , A)$ we have
\begin{align*}
\langle\varphi ,T^{*}f\rangle & = \int_{\Omega}\langle\varphi(\omega),(T^{*}f)(\omega)\rangle d\mu(\omega)\\
& = \int_{\Omega}\langle\varphi(\omega),\langle f, F(\omega)\rangle\rangle d\mu(\omega)\\
& =\int_{\Omega}\varphi(\omega)\langle F(\omega),f\rangle d\mu(\omega)\\
& =\langle \int_{\Omega}\varphi(\omega) F(\omega) d\mu(\omega),f\rangle = \langle T\varphi ,f\rangle.
\end{align*}

(ii)\; The pre-frame operator $T$ is well defined and bounded with $\| T\|\leq\sqrt{B}$ because for $\varphi\in L^{2}(\Omega , A)$ we have
\begin{align*}
\| T\varphi\|^{2} & =\| \int_{\Omega}\varphi(\omega) F(\omega) d\mu(\omega)\|^{2}\\
& =\sup_{f\in U \; ,\; \| f\|=1}\|\langle \int_{\Omega}\varphi(\omega) F(\omega) d\mu(\omega),f\rangle\|\\
& =\sup_{f\in U \; ,\; \| f\|=1}\| \int_{\Omega}\langle\varphi(\omega),\langle f, F(\omega)\rangle\rangle d\mu(\omega)\|\\
& =\sup_{f\in U \; ,\; \| f\|=1}\| \langle\varphi ,T^{*}f\rangle\|\\
&\leq \sup_{f\in U \; ,\; \| f\|=1}\| \langle\varphi ,\varphi\rangle\|\;\|\langle T^{*}f ,T^{*}f\rangle\|\\
& =\sup_{f\in U \; ,\; \| f\|=1}\|\varphi\|^{2}\;\|\int_{\Omega}\langle f,F(\omega)\rangle\langle F(\omega) ,f\rangle d\mu(\omega)\|,
\end{align*}
since $F$ is a continuous frame, so
\begin{equation*}
\|T\varphi\|^{2}\leq \sup_{f\in U \; ,\; \| f\|=1}\|\varphi\|^{2}B\|f\|^{2}\leq B\|\varphi\|^{2},
\end{equation*}
i.e. $\| T\|\leq\sqrt{B}$.

(iii)\; $T$ is surjective. Indeed, by definition of continuous frames in Hilbert $C^{\ast}$-modules, for each $f\in U$,
\begin{equation*}
A\langle f,f\rangle\leq\langle T^{*}f ,T^{*}f\rangle\leq B\langle f,f\rangle .
\end{equation*}
Then $T^{*}$ is bounded below with respect to the inner product and by lemma \ref{SB}, $T$ is surjective.

(iv)\; $T^{*}$ is injective. Indeed, if $f\in U$ and $T^{*}f=0$, then
\begin{align*}
\| \langle f,f\rangle\| & =\| A^{-1}A\langle f,f\rangle\|
 =A^{-1}\| A\langle f,f\rangle\|\\
& \leq A^{-1}\|\langle T^{*}f,T^{*}f\rangle\|= A^{-1}\| T^{*}f\|^{2}.
\end{align*}
Thus $\| \langle f,f\rangle\|=0$ and $f=0$.

(v)\;Now we show that $T^{*}$ has closed range. Let $\lbrace T^{*}f_{n}\rbrace_{n=1}^{\infty}$ be a sequence in $R(T^{*})$ such that $\lim_{n \to \infty} T^{*}f_{n}=g$. By definition of continuous frames in Hilbert $C^{\ast}$-modules, for $n,m\in N$
\begin{equation*}
\Vert A\langle f_{n}-f_{m},f_{n}-f_{m}\rangle\Vert\leq\Vert\langle T^{*}(f_{n}-f_{m}) ,T^{*}(f_{n}-f_{m})\rangle\Vert 
=\Vert T^{*}(f_{n}-f_{m})\Vert^{2}.
\end{equation*}
Since $\lbrace T^{*}f_{n}\rbrace_{n=1}^{\infty}$ is a Cauchy sequence in $L^{2}(\Omega , A)$ so,
\begin{equation*}
\lim_{n,m \to \infty}\Vert A\langle f_{n}-f_{m},f_{n}-f_{m}\rangle\Vert\ =0 .
\end{equation*}
Also
\begin{equation*}
\Vert\langle f_{n}-f_{m},f_{n}-f_{m}\rangle\Vert \leq A^{-1}\Vert A\langle f_{n}-f_{m},f_{n}-f_{m}\rangle\Vert,
\end{equation*}
thus the sequence $\lbrace f_{n}\rbrace_{n=1}^{\infty}$ is a Cauchy sequence in $U$ and so
\begin{equation*}
\exists f\in U\;\;\;s.t\;\;\;\lim_{n \to \infty} f_{n}=f.
\end{equation*}
Definition of continuous frames, implies that
\begin{equation*}
\Vert T^{*}(f_{n}-f)\Vert^{2}\leq B\Vert\langle f_{n}-f ,f_{n}-f\rangle\Vert ,
\end{equation*}
 then $lim_{n \to \infty} \Vert T^{*}f_{n}-T^{*}f\Vert =0$ and $T^{*}f=g$. Therefore $R(T^{*})$ is close.
\end{proof}
\begin{theorem}
Let $F:\Omega \to U$ be a continuous frame for Hilbert $C^{\ast}$-module $U$ with bounds $A,B$ and frame operator $S$ and pre-frame operator $T$. Then  $S=T T^{\ast}$ is positive, adjointable, self-adjoint and invertible and $\| S\|\leq B$.
\end{theorem}

\begin{proof}
Let $F:\Omega \to U$ be a continuous frame for Hilbert $C^{\ast}$-module $U$ with bounds $A,B$. Then

(i)\; $S=T T^{\ast}$ because for $f\in U$ we have
\begin{align*}
\langle T T^{\ast}f,f\rangle & =\langle T(\lbrace\langle f,F(\omega)\rangle\rbrace_{\omega\in\Omega})  ,f\rangle\\
& =\langle\int_{\Omega}\langle f, F(\omega)\rangle F(\omega) d\mu(\omega) ,f\rangle =\langle Sf,f\rangle.
\end{align*}

(ii)\; For $f,g\in U$, we have
\begin{align*}
\langle Sf,g\rangle & =\langle\int_{\Omega}\langle f, F(\omega)\rangle F(\omega) d\mu(\omega) ,g\rangle\\
& =\int_{\Omega}\langle f, F(\omega)\rangle\langle F(\omega),g\rangle d\mu(\omega)\\
& =\int_{\Omega}(\langle g, F(\omega)\rangle\langle F(\omega),f\rangle)^{*} d\mu(\omega)\\
& =(\int_{\Omega}\langle g, F(\omega)\rangle\langle F(\omega),f\rangle d\mu(\omega))^{*}\\
& =(\langle Sg,f\rangle)^{*}=\langle f,Sg\rangle.
\end{align*}
Therefore $S$ is adjointable and self-adjoint.

(iii)\;  By definition of continuous frames, for each $f\in U$ we get 
\begin{equation*}
A\langle f,f\rangle\leq\langle Sf ,f\rangle\leq B\langle f,f\rangle,
\end{equation*}
by \cite[Proposition 2.1.3]{MT}, this implies that $S$ is positive. Also 
\begin{align*}
AI\leq S\leq BI & \Longrightarrow 0\leq I-B^{-1}S\leq I-B^{-1}AI = \dfrac{B-A}{B}I \leq I\\
& \Longrightarrow\Vert I-B^{-1}S\Vert\leq 1.
\end{align*}
Therefore $S$ is invertible.

(iv)\; Since 
\begin{align*}
\| S\| & =\sup_{f\in U \; ,\; \| f\|\leq1}\Vert\langle Sf,f\rangle\Vert\\
& =\sup_{f\in U \; ,\; \| f\|\leq1}\|\int_{\Omega}\langle f, F(\omega)\rangle\langle F(\omega),f\rangle d\mu(\omega)\|\\
& \leq \sup_{f\in U \; ,\; \| f\|\leq1}\| B\langle f,f\rangle\|\leq B.
\end{align*}
Hence, $S$ is bounded.
\end{proof}
\begin{remark}
Let $F:\Omega \to U$ be a continuous frame for Hilbert $C^{\ast}$-module $U$  with the frame operator $S$ and the pre-frame operator $T$. Invertibility of $S$ implies that $T$ is surjective. \\
In fact, since $S$ is invertible so
\begin{equation*}
\forall f\in U\;\;\;\;\;\exists g\in U\;\;\;\;\; s.t\;\;\;\;\; f=Sg=\int_{\Omega}\langle g, F(\omega)\rangle F(\omega) d\mu(\omega) .
\end{equation*}
On the other hand, $S=TT^{*}$ then
\begin{equation*}
 f=TT^{*}g=\int_{\Omega}\langle g, F(\omega)\rangle F(\omega) d\mu(\omega) ,
\end{equation*}
this implies that
\begin{equation*}
 T^{*}g=\lbrace\langle g, F(\omega)\rangle\rbrace_{\omega\in\Omega} \in L^{2}(\Omega , A) ,
\end{equation*}
i.e. $T$ is surjective.
\end{remark}
 \begin{remark}
Let $F:\Omega \to U$ be a continuous frame for Hilbert $C^{\ast}$-module $U$ with the frame operator $S$ and the pre-frame operator $T$. Since $T^{*}:U \to L^{2}(\Omega , A)$ is injective and has closed range so by \cite[Lemma 0.1]{AD}, $S=TT^{*}$ is invertible and\\
$$\Vert S^{-1}\Vert^{-1}I\leq S\leq\Vert T^{*}\Vert^{2}I.$$

Moreover, lower and upper bounds of $F$ are respectively $\Vert S^{-1}\Vert^{-1}$ and $\Vert T\Vert^{2}$.
\end{remark}
 \begin{theorem}
Let $F:\Omega \to U$ be a continuous frame for Hilbert $C^{\ast}$-module $U$ with bounds $A,B$ and the frame operator $S$.
Let $K\in End^{*}(U)$ be invertible.Then $KF$ is a continuous frame for Hilbert $C^{\ast}$-module $U$ with frame operator $KSK^{*}$.
\end{theorem}

\begin{proof}
For each $f\in U$ we have
\begin{align*}
& A\langle K^{*}f,K^{*}f\rangle\leq\int_{\Omega}\langle K^{*}f, F(\omega)\rangle\langle F(\omega),K^{*}f\rangle d\mu(\omega)\leq B\langle K^{*}f,K^{*}f\rangle,
\end{align*}
and so
\begin{align*}
 A\langle K^{*}f,K^{*}f\rangle\leq\int_{\Omega}\langle f, KF(\omega)\rangle\langle KF(\omega),f\rangle d\mu(\omega)\leq B\langle K^{*}f,K^{*}f\rangle
\end{align*}
Since $K$ is invertible so is surjective, then by \cite[Lemma 0.1]{AD},
\begin{align*}
& \Vert K^{-1}\Vert^{-2}=\Vert (KK^{*})^{-1}\Vert^{-1}\leq KK^{*},
\end{align*}
hence,
\begin{align*} \Vert K^{-1}\Vert^{-2}\langle f,f\rangle\leq\langle KK^{*}f,f\rangle =\langle K^{*}f,K^{*}f\rangle\leq\Vert K^{2}\Vert\langle f,f\rangle
\end{align*}
Consequently,
\begin{equation*}
A\Vert K^{-1}\Vert^{-2}\langle f,f\rangle\leq\int_{\Omega}\langle f, KF(\omega)\rangle\langle KF(\omega),f\rangle d\mu(\omega)\leq B\Vert K^{2}\Vert\langle f,f\rangle.
\end{equation*}
Therefore $KF$ is a continuous frame for Hilbert $C^{\ast}$-module $U$.\\
Moreover, for each $f\in U$ we have
\begin{align*}
KSK^{*}f &=K\int_{\Omega}\langle K^{*}f, F(\omega)\rangle F(\omega)d\mu(\omega)\\
& =\int_{\Omega}\langle f,KF(\omega)\rangle KF(\omega)d\mu(\omega),
\end{align*}
thus $KSK^{*}$ is frame operator for $KF$.
\end{proof}
Since $AI\leq S\leq BI$, so  
$$B^{-1}I\leq S^{-1}\leq A^{-1}I.$$
Hence, we get the following corollary.
 \begin{corollary}
Let $F:\Omega \to U$ be a continuous frame for Hilbert $C^{\ast}$-module $U$ with bounds $A,B$ and frame operator $S$. Then $S^{-1}F$ is a continuous frame for Hilbert $C^{\ast}$-module $U$ with bounds $B^{-1},A^{-1}>0$ and frame operator $S^{-1}$.\\
\end{corollary}
\begin{lemma}
Let $U$ be a Hilbert $C^{\ast}$-module. Then the following are equivalent:\\
(i)\;$F:\Omega \to U$ is a continuous Bessel mapping for $U$.\\
(ii)\;The mapping $\Omega \to \langle f,F(\Omega)\rangle$ is measurable and there exists constant $D>0$ such that
\begin{equation*}
\Vert\int_{\Omega}\langle f,F(\omega)\rangle\langle F(\omega),f\rangle d\mu(\omega)\Vert\leq D\Vert \langle f,f\rangle\Vert .
\end{equation*}
\end{lemma}
\begin{proof}
$(i)\Longrightarrow(ii)$ Obvious.\\
$(ii)\Longrightarrow(i)$ Define an operator $V:U \to L^{2}(\Omega , A)$ by $Vf=\lbrace\langle f,F(\Omega)\rangle\rbrace_{\omega\in\Omega}$. It is clear that $V$ is well-defined and $\mathcal A$-linear. Also,
\begin{align*}
\Vert Vf\Vert^{2}=\Vert \langle Vf,Vf\rangle\Vert &=\Vert\int_{\Omega}\langle f,F(\omega)\rangle\langle F(\omega),f\rangle d\mu(\omega)\Vert \\
& \leq D\Vert \langle f,f\rangle\Vert = D\Vert f\Vert^{2}
\end{align*}
Then $V$ is bounded and by theorem \ref{BL} $\langle Vf,Vf\rangle\leq D\langle f,f\rangle$.
\end{proof}
Now, we give an equivalent condition for the mapping $F:\Omega \to U$ to become a continuous frame. 
\begin{theorem}
Let $U$ be a Hilbert $C^{\ast}$-module. Then the following are equivalent:\\
(i)\;$F:\Omega \to U$ is a continuous frame for $U$.\\
(ii)\;The mapping $\Omega \to \langle f,F(\Omega)\rangle$ is measurable and there exist constants $A,B>0$ such that
\begin{equation*}
A\Vert \langle f,f\rangle\Vert\leq\Vert\int_{\Omega}\langle f,F(\omega)\rangle\langle F(\omega),f\rangle d\mu(\omega)\Vert\leq B\Vert \langle f,f\rangle\Vert .
\end{equation*}
\end{theorem}

\begin{proof}
$(i)\Longrightarrow(ii)$ Obvious.\\
$(ii)\Longrightarrow(i)$ Suppose that there exist constants $A,B>0$ such that
\begin{equation*}
A\Vert \langle f,f\rangle\Vert\leq\Vert\int_{\Omega}\langle f,F(\omega)\rangle\langle F(\omega),f\rangle d\mu(\omega)\Vert\leq B\Vert \langle f,f\rangle\Vert .
\end{equation*}
then
\begin{equation*}
\Vert\int_{\Omega}\langle f,F(\omega)\rangle\langle F(\omega),f\rangle d\mu(\omega)\Vert =\Vert \langle Sf,f\rangle\Vert =\Vert \langle S^{\frac{1}{2}}f,S^{\frac{1}{2}}f\rangle\Vert =\Vert S^{\frac{1}{2}}f\Vert^{2}
\end{equation*}
then,
\begin{equation*}
A\Vert f\Vert^{2}\leq\Vert S^{\frac{1}{2}}f\Vert^{2}\;\;\;\Longrightarrow\;\;\;\sqrt{A}\Vert f\Vert\leq\Vert S^{\frac{1}{2}}f\Vert
\end{equation*}
then by Lemma \ref{SB},
\begin{equation*}
\sqrt{A}\langle f,f\rangle\leq\langle S^{\frac{1}{2}}f,S^{\frac{1}{2}}f\rangle =\langle Sf,f\rangle
\end{equation*}
Since $S$ is self-adjoint so by theorem \ref{BL},
\begin{equation*}
\exists D>0\;\;\;\; s.t\;\;\;\; \langle Sf,f\rangle =\langle S^{\frac{1}{2}}f,S^{\frac{1}{2}}f\rangle\leq D\langle f,f\rangle.
\end{equation*}
Hence
\begin{equation*}
\sqrt{A}\langle f,f\rangle\leq\langle Sf,f\rangle =\int_{\Omega}\langle f,F(\omega)\rangle\langle F(\omega),f\rangle d\mu(\omega)\leq D\langle f,f\rangle.
\end{equation*}
\end{proof}
\begin{lemma}
Let $K:U \to L^{2}(\Omega , A)$ be a continuous Bessel mapping for $U$ with bound B and $V:L^{2}(\Omega , A) \to L^{2}(\Omega , A)$ be adjointable. Then $VK:\Omega \to U$ is a continuous Bessel mapping for $U$.
\end{lemma}

\begin{proof}
By Theorem \ref{BL}, we have
\begin{align*}
\int_{\Omega}\langle f,VK(\omega)\rangle\langle VK(\omega),f\rangle d\mu(\omega) & = \int_{\Omega}\langle V^{*}f,K(\omega)\rangle\langle K(\omega),V^{*}f\rangle d\mu(\omega)\\
& \leq B\langle V^{*}f,V^{*}f\rangle\leq BD\langle f,f\rangle,
\end{align*}
for some constant $D>0$.
\end{proof}
\begin{theorem}
Let $U$ be a Hilbert $C^{\ast}$-module over a unital $C^*$-algebra $\mathcal A$. Then a mapping $F:\Omega \to U$ is a continuous frame for $U$ if and only if the synthesis operator $T:L^{2}(\Omega ,A) \to U$ is well-defined and onto.
\end{theorem}

\begin{proof}
$(\Longrightarrow)$ It is shown in Theorem \ref{TH1} .\\
$(\Longleftarrow)$ Let $T$ be well-defined. Then $T$ is adjointable and \\
\begin{align*}
T^{*}: & U\to L^{2}(\Omega ,A)\\
& f\; \longmapsto\;\lbrace\langle f,F(\omega)\rangle\rbrace_{\omega\in\Omega} .
\end{align*}
Also
\begin{align*}
\Vert\int_{\Omega}\langle f,F(\omega)\rangle\langle F(\omega),f\rangle d\mu(\omega)\Vert^{2} & =\Vert\langle f,\int_{\Omega}\langle f,F(\omega)\rangle F(\omega) d\mu(\omega)\rangle\Vert^{2}\\
& \leq \Vert f\Vert^{2} \Vert T(\lbrace\langle f,F(\omega)\rangle\rbrace_{\omega\in\Omega})\Vert^{2}\\
& \leq\Vert f\Vert^{2} \Vert T\Vert^{2} \Vert\lbrace\langle f,F(\omega)\rangle\rbrace_{\omega\in\Omega}\Vert^{2}\\
& = \Vert f\Vert^{2} \Vert T\Vert^{2}\Vert\int_{\Omega}\langle f,F(\omega)\rangle\langle F(\omega),f\rangle d\mu(\omega)\Vert.
\end{align*}
Then
\begin{equation*}
\Vert\int_{\Omega}\langle f,F(\omega)\rangle\langle F(\omega),f\rangle d\mu(\omega)\Vert\leq\Vert f\Vert^{2} \Vert T\Vert^{2} =\Vert\langle f,f\rangle\Vert\Vert T\Vert^{2}.
\end{equation*}
Moreover, since $T$ is onto,so $T^{*}$ is bounded below by lemma \ref{SB}, and $T^{*}\vert_{R(T^{*})}$ is invertible. Then for each $f\in U$ we have $(T^{*}\vert_{R(T^{*})})^{-1}T^{*}f=f$. Then
\begin{align*}
\Vert f\Vert^{2} & = \Vert (T^{*}\vert_{R(T^{*})})^{-1}T^{*}f\Vert^{2}\\
& \leq\Vert (T^{*}\vert_{R(T^{*})})^{-1}\Vert^{2}\Vert T^{*}f\Vert^{2}\\
& = \Vert (T^{*}\vert_{R(T^{*})})^{-1}\Vert^{2}\Vert\int_{\Omega}\langle f,F(\omega)\rangle\langle F(\omega),f\rangle d\mu(\omega)\Vert,
\end{align*}
hence
\begin{equation*}
\Vert (T^{*}\vert_{R(T^{*})})^{-1}\Vert^{-2}\Vert\langle f,f\rangle\Vert\leq\Vert\int_{\Omega}\langle f,F(\omega)\rangle\langle F(\omega),f\rangle d\mu(\omega)\Vert .
\end{equation*}
\end{proof}
 \begin{corollary}
Let $U$ be a Hilbert $C^{\ast}$-module over a unital $C^*$-algebra $\mathcal A$. Then a mapping $F:\Omega \to U$ is a continuous Bessel mapping with bound $B$ for $U$ if and only if the synthesis operator $T:L^{2}(\Omega ,A) \to U$ is well-defined and bounded with $\Vert T\Vert\leq\sqrt{B}$.
\end{corollary}
 
 \section{Duals of continuous frames in Hilbert $C^{\ast}$-modules}
 
 In this section, we introduce the consept of duals of continuous frames in Hilbert $C^{\ast}$-modules and give some important properties of continuous frames and their duals.
\begin{definition}
Let $F:\Omega \to U$ be a continuous Bessel mapping. A continuous Bessel mapping $G:\Omega \to U$ is called a \textit{dual} for $F$ if
\begin{equation*}
f= \int_{\Omega}\langle f,G(\omega)\rangle F(\omega)d\mu(\omega)\;\;\;\;\;(f\in U) ,
\end{equation*}
or
\begin{equation}\label{eq2}
\langle f,g\rangle= \int_{\Omega}\langle f,G(\omega)\rangle\langle F(\omega),g\rangle d\mu(\omega)\;\;\;\;\;(f,g\in U).
\end{equation}
In this case $(F,G)$ is called a \textit{dual pair}. If $T_{F}$ and $T_{G}$ denote the synthesis operators of $F$ and $G$, respectively, then \eqref{eq2} is equivalent to $T_{F}T^{*}_{G}=I_{U}$. The condition 
\begin{equation*}
\langle f,g\rangle= \int_{\Omega}\langle f,G(\omega)\rangle\langle F(\omega),g\rangle d\mu(\omega)\;\;\;\;\;(f,g\in U),
\end{equation*}
is equivalent
\begin{equation*}
\langle f,g\rangle= \int_{\Omega}\langle f,F(\omega)\rangle\langle G(\omega),g\rangle d\mu(\omega)\;\;\;\;\;(f,g\in U),
\end{equation*}
because $T_{F}T^{*}_{G}=I_{U}$ if and only if $T_{G}T^{*}_{F}=I_{U}$ .
\end{definition}
\begin{remark}
Let $F:\Omega \to U$ be a continuous frame for Hilbert $C^{\ast}$-module $U$. Then by reconstructin formula we have
\begin{equation*}
f=S^{-1}Sf=S^{-1}\int_{\Omega}\langle f,F(\omega)\rangle F(\omega) d\mu(\omega)=\int_{\Omega}\langle f,F(\omega)\rangle S^{-1}F(\omega) d\mu(\omega),
\end{equation*}
and
\begin{equation*}
f=SS^{-1}f=S^{-1}\int_{\Omega}\langle S^{-1}f,F(\omega)\rangle F(\omega) d\mu(\omega)=\int_{\Omega}\langle f,S^{-1}F(\omega)\rangle F(\omega) d\mu(\omega).
\end{equation*}
\end{remark}
Then $S^{-1}F$ is a dual for $F$, which is called \textit{canonical dual}.
\begin{definition}
Let $F:\Omega \to U$ be a continuous frame for Hilbert $C^{\ast}$-module $U$. If $F$ has only one dual, we call $F$ a \textit{Riesz-type frame}.
\end{definition}
\begin{example} 
Let $\mathcal A$ and  $(\Omega ,\mu)$ be as in Example \ref{ex1}. Consider the mappings $F,G:\Omega\to \mathcal A$ defined by 
$F(\omega)=\begin{pmatrix}
2\omega&0\\
0&\omega -1
\end{pmatrix}$ and $G(\omega)=\begin{pmatrix}
\dfrac{3}{2}\omega&0\\
0&\omega-\dfrac{7}{3}
\end{pmatrix}$, for any $\omega\in \Omega$.\\ 
Then both $F$ and $G$ are continuous frames for $\mathcal A$ where $\dfrac{1}{3},\dfrac{4}{3}$ are the frame bounds of $F$ and $\dfrac{3}{4},\dfrac{31}{9}$ are the frame bounds of $G$. Also, for each $f=
\begin{pmatrix}
a&0\\
0&b
\end{pmatrix}
\in \mathcal A$, we have
\begin{align*}
\int_{\Omega}\langle f,F(\omega)\rangle G(\omega) d\mu(\omega) &=\int_{[0,1]}
\begin{pmatrix}
2\omega a&0\\
0&(\omega-1) b
\end{pmatrix}\begin{pmatrix}
\dfrac{3}{2}\omega&0\\
0&\omega-\dfrac{7}{3}
\end{pmatrix} d\mu(\omega)\\
&=\begin{pmatrix}
a&0\\
0&b
\end{pmatrix}\int_{[0,1]}
\begin{pmatrix}
3\omega^{2} &0\\
0&\omega^{2}-\dfrac{10}{3}\omega+\dfrac{7}{3}
\end{pmatrix} d\mu(\omega)\\
&=\begin{pmatrix}
a&0\\
0&b
\end{pmatrix}\begin{pmatrix}
1&0\\
0&1
\end{pmatrix} =\begin{pmatrix}
a&0\\
0&b
\end{pmatrix}=f.\\
\end{align*}
Therefore $F$ is a dual of $G$.\\
\end{example}
\begin{theorem}
Let $F:\Omega \to U$ be a continuous frame for Hilbert $C^{\ast}$-module $U$ over a unital $C^*$-algebra $\mathcal A$. Then $F$ is a Riesz-type frame if and only if the analysis operator $T^{*}_{F}:U\to L^{2}(\Omega ,A)$ is onto.
\end{theorem}

\begin{proof}
$(\Longrightarrow)$ Since $T_{F}$ is adjointable and has closed range, so $L^{2}(\Omega ,A)=Ker(T_{F})\oplus R(T^{*}_{F})$ by Lemma \ref{UKR}. Also by \cite[Lemma 15.3.4]{OLS},
\begin{equation*}
R(T^{*}_{F})\neq L^{2}(\Omega ,A)\;\;\;\;\;\Longleftrightarrow\;\;\;\;\;R(T^{*}_{F})^{\perp}\neq\lbrace 0\rbrace
\end{equation*}
Now let $G=S^{-1}F$  be the canonical dual of $F$ and $R(T^{*}_{F})\neq L^{2}(\Omega ,A)$. Then $R(T^{*}_{F})^{\perp}\neq\lbrace 0\rbrace$.\\
Suppose that $h\in R(T^{*}_{F})^{\perp}$ such that $\Vert\int_{\Omega}\vert h(\omega)^{*}\vert^{2} d\mu(\omega)\Vert =1$ i.e. $\Vert h\Vert=1$.\\
Define
\begin{align*}
K: & \Omega\to L^{2}(\Omega ,A)\\
& \omega\; \longmapsto\; h(\omega)^{*}h
\end{align*}
Then for $\varphi\in L^{2}(\Omega ,A)$ we have,
\begin{align*}
\Vert\int_{\Omega}\langle \varphi ,K(\omega)\rangle\langle K(\omega),\varphi\rangle d\mu(\omega)\Vert & =\Vert\int_{\Omega}\langle \varphi ,h(\omega)^{*}h\rangle\langle h(\omega)^{*}h,\varphi\rangle d\mu(\omega)\Vert\\
& =\Vert\int_{\Omega}\langle \varphi ,h\rangle h(\omega)h(\omega)^{*}\langle h,\varphi\rangle d\mu(\omega)\Vert \\
& \leq\Vert\langle \varphi ,h\rangle\Vert\Vert\int_{\Omega}\vert h(\omega)^{*}\vert^{2} d\mu(\omega)\Vert\Vert\langle h,\varphi\rangle\Vert\\
& =\Vert\langle \varphi ,h\rangle\Vert\Vert\langle h,\varphi\rangle\Vert\\
& \leq\Vert h\Vert^{2}\Vert \varphi\Vert^{2} =\Vert h\Vert^{2}\Vert\langle\varphi ,\varphi\rangle\Vert,
\end{align*}
so $K$ is a continuous Bessel mapping.\\
Let $V:L^{2}(\Omega ,A)\to U$ be a adjointable operator such that $Vh\neq 0$. Then $VK:\Omega\to U$ is a continuous Bessel mapping and so $G+VK$ is.\\
If we show that 
\begin{equation*}
\int_{\Omega}\langle f,VK(\omega)\rangle\langle F(\omega),g\rangle d\mu(\omega)=0\;\;\;\;\;\;\;(f,g\in U),
\end{equation*}
then $G+VK$ is a dual of $F$. Because for $f,g\in U$,
\begin{align*}
\int_{\Omega}\langle f,G(\omega)+VK(\omega)\rangle\langle F(\omega),g\rangle d\mu(\omega) & = \int_{\Omega}\langle f,G(\omega)\rangle\langle F(\omega),g\rangle d\mu(\omega) +\int_{\Omega}\langle f,VK(\omega)\rangle\langle F(\omega),g\rangle d\mu(\omega)\\
& = \langle f,g\rangle +0=\langle f,g\rangle.
\end{align*}
i.e. $F$ is not Riesz-type.\\
For this we have,
\begin{align*}
\int_{\Omega}\langle f,VK(\omega)\rangle\langle F(\omega),g\rangle d\mu(\omega) & =\int_{\Omega}\langle f,Vh(\omega)^{*}h\rangle\langle F(\omega),g\rangle d\mu(\omega)\\
& =\int_{\Omega}\langle V^{*}f,h\rangle h(\omega)\langle F(\omega),g\rangle d\mu(\omega)\\
& =\langle V^{*}f,h\rangle\int_{\Omega} h(\omega)\langle F(\omega),g\rangle d\mu(\omega)\\
& =\langle V^{*}f,h\rangle\langle\lbrace h(\omega)\rbrace_{\omega\in\Omega},\lbrace\langle g,F(\omega)\rangle\rbrace_{\omega\in\Omega}\rangle =0.
\end{align*}
Note that $\lbrace h(\omega)\rbrace_{\omega\in\Omega}\in R(T^{*}_{F})^{\perp}$ and $\lbrace\langle g,F(\omega)\rangle\rbrace_{\omega\in\Omega}\in R(T^{*}_{F})$.

$(\Longleftarrow)$ Let $G_{1},G_{2}$ be two duals of $F$ and $G_{1}\neq G_{2}$. Then 
\begin{align*}
\int_{\Omega}\langle f,G_{1}(\omega)-G_{2}(\omega)\rangle\langle F(\omega),g\rangle d\mu(\omega) & = \int_{\Omega}\langle f,G_{1}(\omega)\rangle\langle F(\omega),g\rangle d\mu(\omega) -\int_{\Omega}\langle f,G_{2}(\omega)\rangle\langle F(\omega),g\rangle d\mu(\omega)\\
& = \langle f,g\rangle - \langle f,g\rangle=0.
\end{align*}
Hence
\begin{equation*}
\langle\lbrace\langle f,G_{1}(\omega)-G_{2}(\omega)\rangle\rbrace_{\omega\in\Omega} ,\lbrace\langle g,F(\omega)\rangle\rbrace_{\omega\in\Omega}\rangle =0.
\end{equation*}
Since $\lbrace\langle f,G_{1}(\omega)-G_{2}(\omega)\rangle\rbrace_{\omega\in\Omega}\in R(T^{*}_{G_{1}-G_{2}})$ and $\lbrace\langle g,F(\omega)\rangle\rbrace_{\omega\in\Omega}\in R(T^{*}_{F})$, so
\begin{equation*}
R(T^{*}_{G_{1}-G_{2}})\perp R(T^{*}_{F}).
\end{equation*}
Also $T^{*}_{F}$ is onto and $L^{2}(\Omega ,A)=Ker(T_{F})\oplus R(T^{*}_{F})$, then
\begin{align*}
R(T^{*}_{F})^{\perp}=\lbrace 0\rbrace & \;\Longrightarrow\;\langle f,G_{1}(\omega)-G_{2}(\omega)\rangle =0\;\;(f\in U)\\
& \Longrightarrow\;G_{1}(\omega)-G_{2}(\omega) =0 \;\;(\omega\in\Omega).
& \Longrightarrow\;G_{1}= G_{2}
\end{align*}
Therefore $F$ is Riesz-type.
\end{proof}
We end this section by the following result.    
\begin{corollary}
Let $F:\Omega \to U$ be a Riesz-type frame for Hilbert $C^{\ast}$-module $U$. Then $F(\omega)\neq 0$ for every $\omega\in\Omega$.
\begin{proof}
For this, let $G$ be the canonical dual of $F$ and $\omega_{0}\in\Omega$ such that $F(\omega_{0})=0$ . Define $G_{1}:\Omega \to U$ where $G_{1}(\omega_{0})\neq 0$ and $G_{1}(\omega)=G(\omega)$ for all $\omega\neq\omega_{0}$. Then $G_{1}$ is a continuous Bessel mapping and
\begin{align*}
f & =\int_{\Omega}\langle f,G(\omega)\rangle F(\omega) d\mu(\omega)\\
& =\int_{\lbrace\omega_{0}\rbrace}\langle f,G(\omega)\rangle F(\omega) d\mu(\omega) +\int_{\Omega\setminus\lbrace\omega_{0}\rbrace}\langle f,G(\omega)\rangle F(\omega) d\mu(\omega)\\
& =\int_{\Omega}\langle f,G_{1}(\omega)\rangle F(\omega) d\mu(\omega).
\end{align*}
Hence $G_{1}$ is a dual of $F$ and $F$ is not Riesz-type.
\end{proof}
\end{corollary}

\end{document}